\documentclass{article}
\usepackage{graphicx} 
\usepackage[utf8]{inputenc}
\usepackage{color}
\usepackage{amssymb}
\usepackage{float}
\usepackage{listings}

\usepackage{bm}  
\usepackage{hyperref}       
\usepackage{fullpage}       
\usepackage{url}            
\usepackage{booktabs}       
\usepackage{amsfonts}       
\usepackage{amssymb}        
\usepackage{nicefrac}       
\usepackage{microtype}      
\usepackage{tikz}           
\usepackage{centernot}

\usepackage[textsize=tiny]{todonotes}

\usepackage[normalem]{ulem}
\newcommand{\stkout}[1]{\ifmmode\text{\sout{\ensuremath{#1}}}\else\sout{#1}\fi}

\usepackage{graphicx}
\graphicspath{{rawfigs/}}

\usepackage[ruled,vlined,linesnumbered]{algorithm2e}
\usepackage{amsmath}
\usepackage{amssymb}
\usepackage{amsthm}
\usepackage{bm}
\usepackage{graphicx}
\graphicspath{{figs/}}
\usepackage{mathtools}
\usepackage{subcaption}

\newcommand{\comment}[1]

\theoremstyle{remark}

\newtheorem{theorem}{Theorem}[section]
\newtheorem{lemma}[theorem]{Lemma}
\newtheorem{corollary}{Corollary}

\newtheorem{example}{Example}[section]

\numberwithin{equation}{section}

\renewcommand{\kappa}{\varkappa}
\newcommand{\be}{\begin{equation}}
\newcommand{\ee}{\end{equation}}

\newcommand{\bq}{\begin{eqnarray}}
\newcommand{\eq}{\end{eqnarray}}

\newcommand{\ba}{\begin{array}}
\newcommand{\ea}{\end{array}}

\newcommand{\twotwo}[4]{\left(\begin{array}{cc}#1&#2\\&\\#3&#4\end{array}\right)}

\title{Exact mean and covariance formulas after diagonal transformations of a multivariate normal}
\author{Rebecca Morrison and Estelle Basor}

\begin{document}

\maketitle
\begin{abstract}
 Consider $\bm X \sim \mathcal{N}(\bm 0, \bm \Sigma)$ and $\bm Y = (f_1(X_1), f_2(X_2),\dots, f_d(X_d))$. We call this a diagonal transformation of a multivariate normal.  In this paper we compute exactly the mean vector and covariance matrix of the random vector $\bm Y.$ This is done two different ways: One approach uses a series expansion for the function $f_i$ and the other a transform method. We compute several examples, show how the covariance entries can be estimated, and compare the theoretical results with numerical ones.\\

\end{abstract}
\section{Introduction}
In 1958, Kruskal considered the following problem~\cite{kruskal1958ordinal}: 

\begin{quote}
     Let $\rho_{xy}$ be the correlation of two bivariate normal variables $X$ and $Y$, and now apply a transformation to each variable, such that the resulting marginals are uniform on $[-1,1]$. What is the new correlation between $X'$ and $Y'$, $\rho'_{xy}$?
     \end{quote}

and found---using rather complicated geometric arguments about quadrants---that $\rho'_{xy} = \frac{6}{\pi} \arcsin{(\rho_{xy}/2)}$.
This was again stated in~\cite{clemen1999correlations}, and eventually became a classic result of multivariate statistics. But how can this be proven in a more direct way?

In this paper, we consider a broader problem: Consider $\bm X \sim \mathcal{N}(\bm 0, \bm \Sigma)$ and $\bm Y = (f_1(X_1), f_2(X_2),\dots, f_d(X_d)$. We call this a diagonal transformation of a multivariate normal (this process may also be referred to as a marginal transformation, and $\bm Y$ is sometimes called a nonparanormal distribution~\cite{liu2009nonparanormal}), and then ask: What are the mean vector and covariance matrix of the random vector $\bm Y$?

To fully specify the problem, it remains to specify the set of functions $\{ f_i \}$. In previous work, we computed the mean and covariance exactly for a limited set of functions~\cite{morrison2022diagonal}. There, the functions were assumed to satisfy three constraints: (1)~each transformation function is the same for all variables, i.e., $f = f_i = f_j \,\,\forall i,j$; (2)~$f$ is odd; and (3)~$f$ has uniformly bounded derivatives at 0. In the current work, all of these assumptions are relaxed: the functions can be different and can be odd, even, or a mix. If the transformation function $f_i$ is smooth, then we require $|f_i^{(a)}(0)| \leq C_i\,K_i^a $ for some constants $C_i, K_i$ and for all $a$. But the $f_i$ can also be given as the transform (Fourier or Laplace) of some function in $L^2$. In this way, the theory accommodates discontinuous functions and those with unbounded derivatives. 

Our motivation to compute moments of such distributions originally came from the fields of measure  transport~\cite{baptista2024learning} and learning probabilistic graphical models~\cite{liu2009nonparanormal}. If we write $\bm Y \sim \pi$, then $\pi = D_\sharp L_\sharp \eta$, where $D$ is a (possibly nonlinear) diagonal map of the form $D(x) = (f_1(x_1),\dots,f_d(x_d))$, $L$ is a linear map, and $\eta = \mathcal{N}(\bm 0 , \bm I_d)$. The distribution $L_\sharp \eta = \mathcal{N}(\bm 0, \Sigma)$ is multivariate Gaussian; its graph (of conditional independence) is revealed by the sparsity of $\Sigma^{-1}$ ($X_i$ and $X_j$ are conditionally independent if and only if $\Sigma_{ij}^{-1} = 0$). A diagonal nonlinear map $D$ does not change the graph, but yields a potentially highly non-Gaussian distribution; see Section 3.2 in ~\cite{baptista2024learning}. So distributions of this type are common test cases for non-Gaussian graph learning algorithms. They also provide a type of copula-like description of multivariate distributions: interactions (marginal and conditional independence) are specified through the covariance or the precision (or both~\cite{boege2023geometry}), while marginal behavior is determined with the transformation functions $f_i$.

More broadly, the field of covariance and precision matrix estimation is active and diverse. Recent advancements include high-dimensional covariance estimation with relaxations of sub-Gaussian assumptions~\cite{kuchibhotla2022moving} and for adaptive filtering~\cite{igreja2024analyzing}, optimal covariance estimation with multi-fidelity samples~\cite{maurais2023multi}, differentially private estimators in the small sample regime~\cite{biswas2020coinpress}, optimal precision estimation for compositional data~\cite{zhang2024care}, and fast multipole methods for spatially correlated data assimilation~\cite{hu2024novel}, to name just a few. In general, such estimation methods fundamentally rely on samples; any full-fledged estimator is accompanied by sample size considerations. In contrast, the current paper computes the covariance exactly (without samples), but instead assumes that the Gaussian covariance $\bm \Sigma$ and the set of transformation functions are known. 

The main contribution of this paper gives an explicit series expansion for the covariance between transformed variables $Y_i$ and $Y_j$ given by $\tau_{ij} \coloneqq \mathbb{E}[Y_i Y_j] - \mathbb{E}[Y_i]\mathbb{E}[Y_j]$ in terms of the original covariance between zero-mean Gaussian variables $X_i$ and $X_j$ given by $\sigma_{ij} \coloneqq \mathbb{E}[X_i X_j]$.  For many transformation functions examined here, this series expansion collapses back down to a single closed-form expression; in other cases, we are left with a convergent series such that $\tau_{ij}$ can be approximately computed to high numerical precision with just a few terms. So this theory is particularly useful for applications, including those above, for which we may have good estimates of $\sigma_{ij}$ and $f_i, f_j$, but are potentially faced with limited samples from the joint distribution of $\bm Y$. 

In any case, one may wish to compute the moments of a nonparanormal distribution, just as Kruskal did nearly 70 years ago. We develop general theory for any functions within the given classes, as well as specific results for that classic problem along with many other examples.

The remainder of this paper is organized as follows. Section~\ref{sec:cov} establishes some notation and initial formulas for the transformed moments after a general transformation. This section relies on simply the definition of moments, Gaussian kernels, and some changes in the order of integration. Section~\ref{sec:comp.smooth} shows the series approach for smooth functions and several examples, while section~\ref{sec:comp.trans} shows the transform approach for functions given as Fourier or Laplace transforms. This section also includes functions that are given as Fourier series; the theory here includes a bit of both the series and transform methods. Then, section~\ref{sec:est} provides general bounds of estimates for $\tau_{ij}$ in terms of $\sigma_{ij}$. Section~\ref{sec:gloss} is a glossary of basic properties of the transformed moments, along with a table of some known transformations (for $f_i = f_j$) of common functions. Section~\ref{sec:num} provides several numerical examples, first for many cases in which theoretical and empirical results strongly agree, and also for some cases in which the empirical results suffer due to sample size and the behavior of the particular function. Section~\ref{sec:conc} concludes the paper with some discussion and ideas for extensions.

\section{Moments after general transformation}\label{sec:cov}
Let $\sigma_{ij}$ be the covariance between zero-mean Gaussian variables $X_i$ and $X_j$. Note $\sigma_{ii} =
\sigma_i^2$, where $\sigma_i$ is the standard deviation of $X_i$. We are interested in the
mean vector and covariance matrix of the transformed random vector $Y = (f_1(X_1), f_2(X_2), \dots,
f_d (X_d))$ where we emphasize that each $f_i$ may be a distinct function. 

We assume first that each $f_i$ is such that for any positive $a$, 
\begin{equation}
    \int_{-\infty}^{\infty} |f_i(x)| e^{-a x^2} dx < \infty.
\end{equation}

It is straightforward to verify that 
the mean of the transformed random variables is given by
$$ \nu_i \coloneqq \mathbb{E}[Y_i] = \frac{1}{\sqrt{2\,\pi \,\sigma_{ii}}}\int_{-\infty}^\infty  f_i(x)\, e^{(-x^2/2\sigma_{ii})}\,dx,$$
and that the covariance entries for $i \neq j$ are 
\begin{equation}{\label{tau.int}}\tau_{ij} = c \int_{-\infty}^{\infty}\int_{-\infty}^{\infty} f_i(x)f_j(w) e^{-\frac{\sigma_{jj}x^2}{2 d}}e^{\frac{\sigma_{ij}xw}{ d}}e^{-\frac{\sigma_{ii}w^2}{2 d}}dx\, dw - \nu_i\nu_j,
\end{equation}
where $d = \sigma_{ii}\sigma_{jj} - \sigma_{ij}^2 \,\,\,(d >0) $ and $c = \frac{1}{d^{1/2}2\pi}.$

For $i = j,$ the variance $\tau_{ii}$ of the variable $Y_i$ is 
\begin{equation}{\label{tau.ii}} \frac{1}{\sqrt{2\,\pi \,\sigma_{ii}}}\int_{-\infty}^{\infty} f_i(x)^2 e^{(-x^2/2\sigma_{ii})} dx -\nu_i^2, 
\end{equation}
which is clear from the definition or can be seen from the formula above using a limiting argument making use of the one-dimensional heat kernel. (To see this, let $\sigma_{ij} = \sigma_{ii}$, complete the square in the exponential arguments, and you are left with a heat kernel expression that can be evaluated as $\sigma_{jj}$ approaches $\sigma_{ii}$.)

It is useful to point out that the above formula~(\ref{tau.int}) implies that if 
$f_i$ is an odd function and $f_j$ even then the above integral is zero. This means that 
if we decompose $f_i$ and $f_j$ into odd and even parts as
\[ f_i = f_{io} + f_{ie}, \,\,\,\,f_j = f_{jo} + f_{je},\]
then $\tau_{ij}$ is the same as
\begin{equation}\tau_{ij} = c \int_{-\infty}^{\infty}\int_{-\infty}^{\infty}( f_{io}(x)f_{jo}(w) +  f_{ie}(x)f_{je}(w))e^{-\frac{\sigma_{jj}x^2}{2 d}}e^{\frac{\sigma_{ij}xw}{ d}}e^{-\frac{\sigma_{ii}w^2}{2 d}}dx\, dw - \nu_i\nu_j.
\end{equation}
Thus, the covariance splits nicely into a sum of the covariances from the odd and even parts of the two functions.

\section{Smooth functions: A series approach}\label{sec:comp.smooth}
In this section we compute the covariance for certain smooth functions without making use of (\ref{tau.int}). The main result using this approach is Theorem \ref{main.th} and the formula provided is valid for all $i$ and $j.$

Our goal is to compare the $\tau_{ij}$ with that of $\sigma_{ij}.$ We show how one can write $\tau_{ij}$ as a power series in $\sigma_{ij}.$ In the case of polynomials the series is finite and easy to write down.

Let $f_i$ be a smooth function of $x$ with derivatives at 0 satisfying the condition
\begin{equation}\label{der.est} 
|f_i^{(a)}(0)| \leq C_i\,K_i^a 
\end{equation} 
for some constants $C_i$ and $K_i$ and with Taylor series
\[f_i(x) = \sum_a \frac{f_i^{(a)}(0)}{a!}x^a.\]
We will also repeatedly use the fact that the moments of a zero mean Gaussian random variable $X_i$
with variance $\sigma_{ii}$ are
\begin{align} \label{eq:moms}
    \mathbb{E}_{\rho} &[X_i^k] =\left\{ \begin{array}{ll} 0, & \text{$k$ odd,}\vspace{1em} \\ 
        \sigma_{ii}^{k/2} (k-1)!!, & \text{$k$ even.}
    \end{array} \right.
\end{align}

In the next three subsections, we show how to compute univariate moments, mixed moments, and provide several examples.

\subsection{Univariate moments after transformation}
Directly we have that the first moment of $Y_i$ is
    \begin{align}
      \nu_i \coloneqq \mathbb E_\pi (Y_i) = \mathbb{E}_\rho \left[ f_i (X_i)\right] = \mathbb{E}
        \left[\sum_{k=0}^\infty \frac{f_i^{(k)}(0)}{k!} X_i^k \right]  
        &= \sum_{k=0}^\infty \frac{f_i^{(k)}(0)}{k!} \mathbb{E} \left[X_i^k \right] \nonumber \\
        &= \sum_{\substack{k=0\\\text{even }k}}^\infty \frac{f_i^{(k)}(0) (k-1)!!}{k!}
        \sigma_{ii}^{k/2}.
    \end{align}

    The second moment is 
    \begin{align}
        \tau_{ii} \coloneqq \mathbb E_\pi(Y_i- \nu_i)^2 = \mathbb{E}_\rho \left[ (f_i (X_i) - \nu_i)^2\right] 
        &= \mathbb{E} \left[\sum_{k=0}^\infty \frac{f_i^{(k)}(0)}{k!} X_i^k \sum_{l=0}^\infty
        \frac{f_i^{(l)}(0)}{l!} X_i^l\right] - \nu_i^2\nonumber \\  
        &= \sum_{k=0}^\infty \sum_{l=0}^\infty \frac{f_i^{(k)}(0)f_i^{(l)}(0)}{k!l!} \mathbb{E} \left[X_i^{k+l} \right] - \nu_i^2 \nonumber \\
        &= \sum_{k=0}^\infty \sum_{\substack{l=0\\k+l \text{ even}}}^\infty \frac{f_i^{(k)}(0)
        f_i^{(l)}(0)(k+l-1)!!}{k!l!}
        \sigma_{ii}^{(k+l)/2} - \nu_i^2.
    \end{align}
   
    In the next subsection, we will simplify the above in the general case for mixed moments.

    \subsection{Mixed moments after transformation}
    We compute the mixed moments as
    \begin{align}
    \tau_{ij} \coloneqq (\bm \Sigma_\pi)_{ij} = \mathbb E_\pi \left[(Y_i - \nu_i)(Y_j -
        \nu_j)\right]
        &=\mathbb E_\rho \left[ f_i(X_i)f_j(X_j)\right] - \nu_{ij} \label{eq:mom}\\
        &= \mathbb{E}_\rho \left[ \sum_{k=0}^\infty \sum_{l=0}^\infty
        \frac{f^{i(k)}(0) f^{j(l)}(0)}{k!l!} X_i^k X_j^l  \right]- \nu_{ij}\\
    &= \sum_{k=0}^\infty \sum_{l=0}^\infty \frac{f_i^{(k)}(0)
        f_j^{(l)}(0)}{k!l!} \mathbb{E}_\rho \left[ X_i^k X^l_j  \right]- \nu_{ij}.
    \end{align}
    where we also set $\nu_{ij} \coloneqq \nu_i \nu_j$.

Let's check that the series above converges.
\begin{lemma}
    Suppose that the derivatives of each $f_i$ satisfy (\ref{der.est}). Then the series
    \[\sum_{k=0}^\infty \sum_{l=0}^\infty \frac{f_i^{(k)}(0)
        f_j^{(l)}(0)}{k!l!} \mathbb{E}_\rho \left[ X_i^k X^l_j  \right]\]
        converges.
\end{lemma}
\proof
Consider
    \begin{equation}\sum_{k=0}^\infty \sum_{l=0}^\infty \frac{f_i^{(k)}(0)
    f_j^{(l)}(0)}{k!l!} \mathbb{E}_\rho \left[ X_i^k X^l_j  \right]\label{eq:mixed}\end{equation}
for fixed $i$ and $j$, and first write 
\[\mathbb{E}_\rho \left[X_i^k X_j^l \right] = \int X_i^k X_j^l \rho(X) dX .\]
By the Cauchy-Schwarz inequality, this square of this expectation is bounded by the product of
integrals:
\begin{align*} \left(\int X_i^k X_j^l \rho(X) dX  \right)^2  &\leq \left(\int X_i^{2k} \rho(X)dX \right)
\left(\int X_j^{2l} \rho(X)dX \right)\\
&\leq \sigma_{ii}^k (2k - 1)!! \sigma_{jj}^l (2l-1)!!.
\end{align*}
Next, let $M$ be a bound on $\sigma_{ii}$ for all $i$, and suppose $|f_i^{(k)}(0)| \leq C\,N^k$ for all $i$
and $k$. Then we have that the sum (\ref{eq:mixed}) is bounded by a constant times
\[\sum_{k=0} \frac{(MN)^k ((2k-1)!!)^{1/2}}{k!}\sum_{l=0} \frac{(MN)^l ((2l-1)!!)^{1/2}}{l!}\]
Note $(2k-1)!! = \frac{(2k-1)!}{(k-1)!2^{k-1}} = \frac{(2k)!}{k! 2^k}$.
We can take one of the sums above:
\[\sum_{k=0} \frac{(MN)^k}{k!}\left( \frac{(2k)!}{k!2^k} \right)^{1/2}\]
which easily converges by the ratio test. Thus the product of both sums converges, and the proof is
complete.
\endproof

Next, based on Wick's theorem~\cite{wick1950evaluation}, or Isserlis's counting theorem~\cite{isserlis1918formula}, we also obtain an
explicit formula for the second-order mixed moments of Gaussian random variables:
\begin{align} \label{eq:wicks}
    \mathbb{E}_{\rho} &[X_i^p X_j^q] =\left\{ \begin{array}{ll} 0, & \text{$p+q$ odd,}\vspace{1em} \\ 
    \sum_{\substack{k = p\\ \text{ by }-2}}^0 (p-k-1)!! {p \choose k}{q \choose k} k! (q-k-1)!!\, \sigma_{ii}^{(p-k)/2} \sigma_{ij}^k
        \sigma_{jj}^{(q-k)/2}, & \text{$p+q$ even.}
    \end{array} \right.
\end{align}

As a final piece of the set up, we introduce a new function
\[F_{ki}(x) = \sum_{u\geq 0} \frac{f_i^{(2u + k)}(0)}{u!} x^u\]
and also $G_{kij}(x) = F_{ki}(\sigma_{ii}x)F_{kj}(\sigma_{jj}x)$.

\begin{lemma}\label{lem:zeroth}
    The mean of the transformed variable $Y_i$ is $\nu_i = F_{0i}(\sigma_{ii}/2)$. 
\end{lemma}
\proof
\begin{align*}
    \nu_i &= \sum_{k \text{ even}} f^{(k)}(0) \sigma_{ii}^{k/2} (k-1)!!\\
    &= \sum_s \frac{f^{(2s)}(0)}{(2s)!} \sigma_{ii}^s (2s -1)!!\\
    &=\sum_s \frac{f^{(2s)}(0)}{s!}\left(\frac{\sigma_{ii}}{2}\right)^s\\
    &= F_{0i}(\sigma_{ii}/2).
\end{align*}
\endproof
\begin{corollary}
From this it immediately follows that $\nu_{ij} = \nu_i \nu_j = G_{0ij}(1/2)$.
\end{corollary}

%
The theorem below yields the power series expansion for $\tau_{ij}.$ The proof of the theorem shows how to explicity construct the power series. 

\begin{theorem}{\label{main.th}} Let $f_i$ be a function with derivatives satisfying (\ref{der.est}), with Taylor series
\[f_i(x) = \sum_a \frac{f_i^{(a)}(0)}{a!}x^a.\]
    Define $F_{ki}$ and $G_{kij}$ the same as above, that is, as
\[F_{ki}(x) = \sum_{u\geq 0} \frac{f_i^{(2u + k)}(0)}{u!} x^u\]
and $G_{kij}(x) = F_{ki}(\sigma_{ii}x)F_{kj}(\sigma_{jj}x)$.
    Then $\sigma_{ij}$ is transformed to \[\tau_{ij} = \sum_{k \geq 1} G_{kij}(1/2)
    \frac{\sigma_{ij}^k}{k!}.\]
\end{theorem}
\proof
\begin{align*}
    \mathbb{E}_\rho [f_i(X_i)f_j(X_j)] + \nu_{ij} &= \mathbb{E} \left[ \sum_{p\geq0} \sum_{q \geq
    0}\frac{f_i^{(p)}(0) f_j^{(q)}(0)}{p! q!} X_i^p X_j^q \right]\\
&= \sum_{p\geq0} \sum_{q \geq
    0}\frac{f_i^{(p)}(0) f_j^{(q)}(0)}{p! q!} \mathbb{E} \left[ X_i^p X_j^q \right]
    \end{align*}
    Now we use the result by Wick's theorem, Eq.~\ref{eq:wicks}:
\begin{align*}
    \mathbb{E}_\rho [f_i(X_i)f_j(X_j)] 
    &= \sum_{p\geq0} \sum_{\substack{q \geq
    0\\p+q \text{ even}}} \frac{f_i^{(p)}(0) f_j^{(q)}(0)}{p! q!} \sum_{\substack{k = p\\ \text{ by }-2}}^0 (p-k-1)!! {p \choose k}{q \choose k} k! (q-k-1)!!\, \sigma_{ii}^{(p-k)/2} \sigma_{ij}^k
        \sigma_{jj}^{(q-k)/2}\\
        &= \sum_{k=0}^\infty \frac{\sigma_{ij}^k}{k!} \sum_{p \geq k} \sum_{\substack{q \geq k\\p+q
        \text{ even}}} \frac{f_i^{(p)}(0)f_j^{(q)}(0) \sigma_{ii}^{(p-k)/2}
        \sigma_{jj}^{(q-k)/2}}{k!2^{(p+q)/2-k}((p-k)/2)!((q-k)/2)!}
    \end{align*}
Reindex by $p -k = 2m$, $q-k = 2n$, then the sums separate and the $k$th coefficient is:
\begin{equation*}
    \sum_{m\geq 0} \frac{f_i^{(2m+k)}(0)}{2^m m!} \sigma_{ii}^m \sum_{n\geq 0}
    \frac{f_j^{(2n+k)}(0)}{2^n n!} \sigma_{jj}^n.
\end{equation*}
Recall $F_{ki}(\sigma_{ii}x) = \sum_{u=0}^\infty \frac{f_i^{(2u+k)}}{u!} x^u \sigma_{ii}^u$ so that
$F_{ki}^{(l)}(0) = f_i^{(2l+ k)}(0)\sigma_{ii}^l$, and so the coefficient is
\[\sum_m \frac{F_{ki}^{(m)}(0)}{2^m m!}\sigma_{ii}^m\sum_n \frac{F_{kj}^{(n)}(0)}{2^n
n!}\sigma_{ii}^n = F_{ki}(\sigma_{ii}/2)
F_{kj}(\sigma_{jj}/2) = G_{kij}(1/2).\]
At this point, we have
\[\tau_{ij} + \nu_{ij} = \sum_{k=0}^\infty G_{kij}(1/2) \frac{\sigma_{ij}^k}{k!}.\]
But, by Lemma~\ref{lem:zeroth}, $\nu_{ij}$ is the 0th order term in the sum, so
\[\tau_{ij} = \sum_{k=1}^\infty G_{kij}(1/2) \frac{\sigma_{ij}^k}{k!}.\]
\endproof

\subsection{Examples, series approach}\label{ssec:ex-d}
This first example compares the two computations of the transformed mean when $f_i = \cos x$.
\begin{example}[Mean, $f_i = \cos x$.]\label{ex:cos-mean}
    First, we can compute the mean term directly:
    \begin{align*} \nu_i &= \sum_{k=0}^\infty \frac{f^{(k)}}{k!} \mathbb{E}_\rho [X_i^k] \\
        &= \sum_{k \text{ even}} \frac{f^{(k)}}{k!} \sigma_{ii}^{k/2} (k-1)!! \\
        &= \sum_l \frac{f^{(2l)}}{(2l)!} \sigma_{ii}^{l} (2l-1)!! \\
        &= \sum_{l=0} \frac{(-\sigma_{ii}/2)^l}{l!} \\
        &= e^{-\sigma_{ii}/2}.
    \end{align*}
    But we could also simply compute this term using Lemma~\ref{lem:zeroth}:
    \begin{align*}
        F_{0i}(x) &= \sum_{u \geq 0} \frac{f_i^{(2u)}(0)x^u}{u!}\\
        &=\sum_{u \geq 0} \frac{(-x)^u }{u!} = e^{-x}.
    \end{align*}
    So $\nu_i = F_{0i}(\sigma_{ii}/2) = e^{-\sigma_{ii}/2}$.
\end{example}

For this same function let us now compute $\tau_{ij}$.
\begin{example}[Covariance, $f_i = f_j = \cos x$.]\label{ex:cos}
    Here $F_{ki} (x)  = (-1)^{k/2}e^x$ if $k$ is even, and 0 otherwise, so $G_{kij}(1/2) =
    e^{-(\sigma_{ii}+\sigma_{jj})/2}$ if $k$ is even, and 0 otherwise. Thus
    \begin{align*} \tau_{ij} &=  \sum_{\substack{k\geq1\\k \text{ even}}} G_{kij}(1/2)
    \frac{\sigma_{ij}^k}{k!}\\
        &= e^{-(\sigma_{ii} + \sigma_{jj})/2}\sum_{\substack{k\geq1\\k \text{ even}}}
        \frac{\sigma_{ij}^k}{k!}\\
        &=e^{-(\sigma_{ii} + \sigma_{jj})/2} (\cosh \sigma_{ij} - 1).
    \end{align*}
\end{example}

\begin{example}[$f = x^3 + x^2 + x$] We again compute the covariance entries $\tau_{ij}$ both ways, first directly (without the use of the previous theorem) and second with Theorem~\ref{main.th}.

    In this example, $f^{(0)}(0) = 0$, $f^{(1)}(0) = 1$, $f^{(2)}(0) = 2$, $f^{(3)}(0) = 6$, and all
    higher derivatives are 0.
    Let's compute the mean term first.
    \begin{align*}
        \nu_i = \sum_k \frac{f^{(k)}(0)}{k!} \mathbb{E}_\rho [X_i^k] &= \sum_{k \text{ even}}
        \frac{f^{(k)}(0)}{k!} \sigma_{ii}^{k/2} (k-1)!!\\
        &= 0 + \frac{2}{2!}\sigma_{ii} + 0\dots\\
        &= \sigma_{ii}.
    \end{align*}

    Next, the diagonal term of the new covariance, $\tau_{ii}$, is:
\begin{align*}
    \tau_{ii} + \nu_{ii} &= \sum_k \sum_l \frac{f^{(k)}(0)f^{(l)}(0)}{k!l!} \mathbb{E}_\rho
    [X_i^{k+l}], \quad \text{$k + l$ even}\\
    &= \mathbb{E}_\rho[X_i^2] + \frac{2 \cdot 6}{1! 3!}\mathbb{E}_\rho [X_i^4] + \frac{2 \cdot 2}{2!
    2!} \mathbb{E}_\rho [X_i^4] + \frac{2 \cdot 6^2}{3! 3!} \mathbb{E}_\rho [X_i^6]\\
    &= \sigma_{ii} + 9 \sigma_{ii}^2 + 15\sigma_{ii}^3.
\end{align*}
    Subtracting off the mean-squared term, we have
    \[ \tau_{ii} = \sigma_{ii} + 8 \sigma_{ii}^2 + 15\sigma_{ii}^3.\]
    Next, we find the off-diagonals of $\Sigma_\pi$, $\tau_{ij}$:
\begin{align*}
    \tau_{ij} + \nu_{ij} &= \sum_k \frac{\sigma_{ij}^k}{k!} \sum_{n=0} \frac{f^{(2n + p)}(0)}{2^n
    n!} \sigma_{ii}^n \sum_{m=0} \frac{f^{(2m + p)}(0)}{2^m
    m!} \sigma_{jj}^m \\
    &= 1 *\left(\left(\frac{f^{(2)}(0)}{2}\right)^2 \sigma_{ii} \sigma_{jj} \right) 
    + \sigma_{ij}^1 *\left(f^{(1)}(0) \frac{f^{(3)}(0)}{2} \sigma_{jj}^1  +  \frac{f^{(3)}(0)}{2} f^{(3)}(0)
    \sigma_{ii}^1  + \left(f^{(3)(0)}{2^2}\sigma_{ii}\sigma_{jj} \right)^2 \right) + \dots\\
    &\quad+\frac{\sigma_{ij}^2}{2} * \left( (f^{(2)} (0))^2 \right) + \frac{\sigma_{ij}^3}{6}
    \left((f^{(3)})^2   \right)\\
    &= \sigma_{ii}\sigma_{jj} + \sigma_{ij} + 3 \sigma_{ii}\sigma_{ij} + 3 \sigma_{jj}\sigma_{ij} +
    2 \sigma_{ij}^2 + 9 \sigma_{ii}\sigma_{jj}\sigma_{ij} + 6 \sigma_{ij}^3 
\end{align*}
and so,
    \[ \tau_{ij} = \sigma_{ij} + 3 \sigma_{ii}\sigma_{ij} + 3 \sigma_{jj}\sigma_{ij} +
    2 \sigma_{ij}^2 + 9 \sigma_{ii}\sigma_{jj}\sigma_{ij} + 6 \sigma_{ij}^3 \]
where, in the last line, we used $\nu_{ij} = \sigma_{ii}\sigma_{jj}$. Also, setting $i=j$ recovers
the result for the diagonal entries.

    On the other hand, we can compute all of the above from Lemma~\ref{lem:zeroth} and Theorem~\ref{main.th}.
    Here $F_{0i}(x) = 2x$, and so $\nu_i = 2{\sigma_{ii}/2} = \sigma_{ii}$.

    For the mixed moment, we have
    \[ F_{1i}(x) = 1 + 6x,\quad F_{2i} (x) = 2, \quad F_{3i} (x) = 6\]
    and thus
    \[ G_{1ij}(1/2) = (1+3\sigma_{ii})(1 + 3\sigma_{jj}), \quad G_{2ij}(1/2) = 4, \quad G_{3ij}(1/2)
    = 36\]
    and from this, 
    \begin{align*}
        \tau_{ij} &=\sum_{k\geq1} G_{kij}(1/2) \frac{\sigma_{ij}^k}{k!}\\
        &= (1+ 3\sigma_{ii})(1 + 3\sigma_{jj}) \sigma_{ij} + 4 \frac{\sigma_{ij}^2}{2!} + 36
        \frac{\sigma_{ij}^3}{3!}\\
        &= (1+  3\sigma_{ii} + 3\sigma_{jj} + 9\sigma_{ii}\sigma_{jj})\sigma_{ij} +2\sigma_{ij}^2 +
        6\sigma_{ij}^3
    \end{align*}
\end{example}
which matches what we found through direct computation above.


\section{Other descriptions of covariance entries via transforms}\label{sec:comp.trans}
In this section we show in three different ways how the transformed $\sigma_{ij}$ is explicitly written in terms of original transformation functions. We make use of the Fourier transform, the Fourier series and the Laplace transform of the coordinate-wise transformation functions.

 \subsection{The Fourier transform}
Let $$f_j(x) = \hat{g}_j(x) = \frac{1}{2\pi}\int_{-\infty}^\infty g_j(y) e^{-i\,x\,y}\,dy,$$ where $g$ is assumed to be in $L^2(-\infty, \infty).$ 

We recall that the inverse transform is given by the formula
$$g_j(x) = \check{f}_j(x) = \int_{-\infty}^\infty f_j(y) e^{i\,x\,y}\,dy,$$
and for later use note that 
$$\hat{\hat{g}}(x) = \frac{1}{2\pi}\,g(-x).$$

Substituting $$f_i(x) = \frac{1}{2\pi}\int_{-\infty}^\infty g_i(y) e^{-i\,x\,y}\,dy, \quad 
f_j(w) = \frac{1}{2\pi}\int_{-\infty}^\infty g_j(z) e^{-i\,w\,z}\,dz,
$$ in (\ref{tau.int}) and changing the order of integration yields
\begin{equation}{\label{tau.fint}}\tau_{ij} = \frac{1}{4\pi^2} \int_{-\infty}^\infty \int_{-\infty}^\infty g_i(y)\,g_j(z)\,(e^{- y\,z\,\sigma_{ij}}-1)\, e^{(-y^2\sigma_{ii} -z^2 \sigma_{jj})/2}\,dydz.
\end{equation}

Let us now show how (\ref{tau.fint}) can also yield a nice series expansion for $\tau_{ij}.$ 
\begin{theorem}\label{thm.trans} Let $f_i$ and  $f_j$ be in $L^2(-\infty, \infty)$ with inverse transform $g_i$ and $g_j$ and with $i \neq j.$ Suppose both $g_i$ and $g_j$ are bounded by a constant $N.$ 
Define
\[F^*_{ki}(x) = \frac{1}{2\pi}(-i)^k \int_{-\infty}^\infty g_i(y)\,\,y^{k} e^{-y^2 x} dy\] and
\[G^*_{kij} = F^*_{ki}(\sigma_{ii}x)F^*_{kj}(\sigma_{jj}x)\] Then

\[\tau_{ij} = \sum_{k=1}^{\infty} \frac{ G^*_{kij}(1/2)(\sigma_{ij})^{k}}{k!}.\]
\end{theorem}
\proof
Expanding the middle exponential in (\ref{tau.fint})yields
$$\tau_{ij}  = \frac{1}{4\pi^2}\sum_{k\geq 1} \frac{(-\sigma_{ij})^k}{k!}\int_{-\infty}^\infty \int_{-\infty}^\infty g_i(y)\,g_j(z)\,y^k z^k\, e^{(-y^2\sigma_{ii} -z^2 \sigma_{jj})/2}\,dydz.$$ We need to show convergence of this sum.
The integral
\[\int_{-\infty}^\infty  g_i(y)\,\,y^k \, e^{(-y^2\sigma_{ii})/2}\,dy\]
is at most 
\[2\,N\int_{0}^\infty\,\,y^k \, e^{(-y^2\sigma_{ii})/2}\,dy\] If $k$ is even this is 
\[2\,N (k-1)!!\sigma_{ii}^{(-(k+1)/2)}.\] If $k$ is odd this is at most
\[2 \, N 2^{(k-1)/2} ((k-1)/2)!\sigma_{ii}^{(-(k+1)/2)}\] With these estimates we have that convergence holds as long as $\frac{\sigma_{ij}^2}{\sigma_{ii}\sigma_{jj}} < 1.$
\endproof
Note that when $i = j$, $\tau_{ii}$ can be computed either with (\ref{tau.ii}) or (\ref{tau.fint}), and in some cases, the series expansion above may also hold.

To show that this agrees with our previous result assuming the Fourier transform representation as above and assuming that our $g_i$s have compact support, we note that 
$$f_{j}^{(2u+k)}(x) = \frac{(-i)^k}{2\pi}\int_{-\infty}^\infty g_j(y)(-1)^u y^{2u+k} e^{-i\,x\,y} \,dy$$
and 
$$f_j^{(2u+k)}(0) = \frac{(-i)^k}{2\pi}\int_{-\infty}^\infty g_j(y)(-1)^u y^{2u+k}\,dy. $$

Thus using our recipe we find that 
$$F_{kj}(x) = \frac{(-i)^k}{2\pi} \int_{-\infty}^\infty g_j(y)\,y^{k}\, e^{-y^2\,x}\,dy $$
and
$$G_{kij}(x) = \frac{(-1)^k}{4\pi^2} \int_{-\infty}^\infty \int_{-\infty}^\infty g_i(y)\,g_j(z)\,y^{k}z^{k}\, e^{(-y^2\sigma_{ii} -z^2 \sigma_{jj})x}\,dydz. $$

Thus we have that $F^*_{kj} = F_{kj},\,\,\,\, G^*_{kij} = G_{kij}$ and  
$$\tau_{ij} = \frac{1}{4\pi^2} \int_{-\infty}^\infty \int_{-\infty}^\infty g_i(y)\,g_j(z)\,(\exp(- y\,z\,\sigma_{ij})-1)\, e^{(-y^2\sigma_{ii} -z^2 \sigma_{jj})/2}\,dydz.$$

The careful reader may wonder why the condition of compact support was required instead of the possibly weaker condition of (\ref{der.est})
$$\Big|\frac{(-i)^k}{2\pi}\int_{-\infty}^\infty g_l(y) y^{k}\,dy\Big|\,\,\leq C_l M_l^k .$$  This is because requiring the above estimate implies that $f_l$ is a entire function of exponential type and thus has an inverse transform with compact support and thus the conditions are equivalent.
However, the formula (\ref{tau.fint}) easily extends to distributions with compact support, for example, any distribution that is the sum of Dirac deltas. 

\subsubsection{Examples, Fourier transform}\label{sssec:ex-ft}
Here is a simple example that uses the formula~\ref{tau.fint}.
\begin{example}[Gaussian function]\label{ex:gauss}
Let $g_i(x) = e^{-\frac{a x^2}{2}}$ and $g_j(x) = e^{-\frac{b x^2}{2}}$ with $a, b$ positive. Here 
\[ f_i(x) = \frac{1}{\sqrt{2\pi\, a}} e^{-\frac{x^2}{2a}},\,\,\,\,f_j(x) = \frac{1}{\sqrt{2\pi\,b}} e^{-\frac{x^2}{2b}}.\]We note that this pair of functions does not satisfy (\ref{der.est}) and thus the method of the previous section does not apply. 

But with our integral formula we have
\[ \tau_{ij} = \frac{1}{4\pi^2} \int_{-\infty}^\infty \int_{-\infty}^\infty \,e^{-\frac{a y^2}{2}}e^{-\frac{b z^2}{2}}(e^{- y\,z\,\sigma_{ij}}-1)\, e^{(-y^2\sigma_{ii} -z^2 \sigma_{jj})/2}\,dydz\]or
\[ \tau_{ij} = \frac{1}{4\pi^2} \int_{-\infty}^\infty \int_{-\infty}^\infty \,(e^{- y\,z\,\sigma_{ij}}-1)\, e^{(-(a +\sigma_{ii})y^2 - (b+ \sigma_{jj})z^2/2}\,dydz\]
which is the same as
\[ \frac{1}{4\pi^2} \int_{-\infty}^\infty \int_{-\infty}^\infty \,(e^{(- y\,z\,\sigma_{ij})}\, e^{(-(a +\sigma_{ii})y^2 - (b+ \sigma_{jj})z^2)/2}\,dydz - \frac{1}{4\pi^2} \int_{-\infty}^\infty \int_{-\infty}^\infty \,\, e^{(-(a +\sigma_{ii})y^2 - (b+ \sigma_{jj})z^2)/2}\,dydz.\]
But the integral on the left corresponds to integrating
$\exp{( -\frac{1}{2}x^T A x)}$ with \[ A = \twotwo{a+\sigma_{ii}}{-\sigma_{ij}}{-\sigma_{ij}}{b+\sigma_{jj}}\]and thus the integral is 
\[\frac{1}{2\pi}((a+\sigma_{ii})(b+\sigma_{jj}) - \sigma_{ij}^2)^{-1/2}.\] One can compute the integral on the right in the same way to find that
\[ \tau_{ij} = \frac{1}{2\pi}(((a+\sigma_{ii})(b+\sigma_{jj}) - \sigma_{ij}^2)^{-1/2} - ((a+\sigma_{ii})(b+\sigma_{jj}))^{-1/2}).\]
\end{example}

Next is an example that illustrates the series in theorem~\ref{thm.trans}. The transform takes continuous data and maps it to binary values, as might occur during classification or other discrete tasks.

\begin{example}[Charateristic function]\label{ex:char}
Let $f = f_i = f_j = \chi_{[-1,1]},$ that is, the function that is $1$ for values between $1$ and $-1$ and zero otherwise. Suppose also that $\sigma_{ii} = \sigma_{jj} = 1.$ In this case  $g(y) = \frac{2\, \sin y}{y}.$ So we have that 
\[ \tau_{ij} = \frac{1}{\pi^2} \int_{-\infty}^\infty \int_{-\infty}^{\infty} \frac{\sin y}{y}\frac{\sin z}{z}\,(\exp(- y\,z\,\sigma_{ij})-1)\, e^{-(y^2 +z^2)/2}\,dydz.\]
Expanding the $(\exp(- y\,z\,\sigma_{ij})-1)$ term and using the fact that $\sin y/ y$ is an even function we have that this is the same as
$$\frac{1}{\pi^2} \sum_{k = 1}^{\infty}\frac{\sigma_{ij}^{2k}}{2k!}\left(\int_{-\infty}^{\infty}\sin y \, y^{2k-1} e^{-y^2/2} dy\right)^2.$$
Now 
\[F_{2k,j}(1/2) = \frac{(-1)^k}{2 \pi} \int_{-\infty}^{\infty}\sin y \,\, y^{2k-1} \,e^{-y^2/2} \,dy\] and the integral above is known to be $\sqrt{\pi}\,H_{2k-1}(\frac{1}{\sqrt{2}})/(\sqrt{\,e}\,\,2^{k-1}),$
where $H_n$ is the Hermite polynomial of order degree $n$ and normalized with leading coefficient $2^n.$ Thus our final answer is
$$\frac{1}{\pi \, e}\sum_{k = 1}^{\infty}\,\frac{\sigma_{ij}^{2k}}{2k!\,\,2^{2k-2}}\,(H_{2k-1}(1/\sqrt{2}))^2.$$
The first two terms of this are
\[ \frac{1}{\,\pi \, e}\left(\frac{\sigma_{ij}^{2}}{2!}(H_{1}(1/\sqrt{2}))^2 + \frac{\sigma_{ij}^{4}}{4\,4!}(H_{3}(1/\sqrt{2}))^2  + \cdots \right) =
 \frac{1}{\,\pi \, e}(\sigma_{ij}^{2} + \frac{\sigma_{ij}^{4}}{3} + \cdots ).\]
The series above converges, as we know, for $|\sigma_{ij}| < 1,$ but this can also be seen by using the asymptotic expansion for the Hermite polynomials. The series is easy to compute using something like mathematica. 
For $i = j$ one can use formula (\ref{tau.ii}) and find that the answer to be $\Phi(1/\sqrt{2}) - \Phi(1/\sqrt{2})^2$ where $\Phi(z) = \frac{2}{\sqrt{\pi}}\int_0^z e^{-t^2} dt,$ (the error function) which is around $0.22$.
\end{example}

\begin{example}[Identity on an interval] \label{ex:char-x}
$f = f_i = f_j = x$ if $-1 < x <1$ and $0$ otherwise and assume that $\sigma_{ii} = \sigma_{jj} = 1.$ Then 
$$g(y) = 2i\left(\frac{\sin y - y \cos y}{y^2}\right).$$ Since $g$ is odd and using the same approach as above we must compute
\[\int_{-\infty}^{\infty}\frac{\sin y - y \cos y}{y^2}  \,\, y^{2k+1} \,e^{-y^2/2} \,dy.\]
If $k = 0$, then we have 
\[\int_{-\infty}^{\infty}\frac{\sin y - y \cos y}{y}  \,\, \,e^{-y^2/2} \,dy,\] or
by direct computation
$$\pi\Phi(1/\sqrt{2}) - (\frac{2\pi}{e})^{1/2},$$ where
$\Phi(z) = \frac{2}{\sqrt{\pi}}\int_0^z e^{-t^2} dt,$ is the error function.
For $k >1$ the integral is
\[\int_{-\infty}^{\infty}(\sin y - y \cos y )  \,\, y^{2k-1} \,e^{-y^2/2} \,dy\]
\[ = \int_{-\infty}^{\infty}\sin y \,\, y^{2k-1} \,e^{-y^2/2}\,dy -  \int_{-\infty}^{\infty}\cos y   \,\, y^{2k} \,e^{-y^2/2} \,dy.\] Both of the last two integrals can be written in terms of Hermite polynomials to find that the above is
\begin{equation}
    (-1)^{k+1}2^{-k+1}\sqrt{\frac{2\pi}{e}}(2k-1)H_{2k-2}(1/\sqrt{2}).
\end{equation}
This gives a series of 
$$\frac{1}{\pi^2}\left(\pi\Phi(1/\sqrt{2}) - \left(\frac{2\pi}{e}\right)^{1/2}\right)^2\sigma_{ij}+ \frac{2}{ \pi e} \sum_{k=2}^{\infty} \frac{\sigma_{ij}^{2k+1}(2k-1)^2}{(2k+1)!2^{2k-2}}H_{2k-2}(1/\sqrt{2})^2.$$
For $i = j$, (\ref{tau.ii}) gives $\tau_{ii} = \Phi(1/\sqrt{2}) - \sqrt{\frac{2}{\pi e}} \approx 0.20$. 
\end{example}

\begin{example}\label{ex:uni}[Normal to uniform]
Here we compute the covariance for the case when $\sigma_{ii} = \sigma_{jj} = 1$ and
\[f_i(x) = f_j(x) = f(x) = \frac{1}{\sqrt{2\pi}} \int_{-\infty}^x e^{-t^2/2}.\]
This is the cumulative distribution function of a univariate normal, so this maps each marginal to uniform on $[0,1]$. In the principle-valued sense the corresponding $g$ is given by
$ \frac{e^{-y^2/2}}{-iy} - \frac{1}{2}\delta_0.$ The delta part does not contribute to the computation and thus 
\[\tau_{ij} = \frac{1}{4\pi^2} \int_{-\infty}^\infty \int_{-\infty}^\infty \frac{e^{-y^2/2}}{iy} \frac{e^{-z^2/2}}{iz}\,(e^{- y\,z\,\sigma_{ij}}-1)\, e^{(-y^2\ -z^2 )/2}\,dydz.\] 
If we think of $\tau_{ij}$ as a function of $\sigma_{ij}$ and differentiate we find that 
the derivative is given by
\[\frac{1}{4\pi^2} \int_{-\infty}^\infty \int_{-\infty}^\infty e^{-y^2/2} e^{-z^2/2}\,e^{- y\,z\,\sigma_{ij}}\, e^{(-y^2 -z^2 )/2}\,dydz.\]
From our previous computation we know this is 
\[\frac{1}{2\pi} (4 - \sigma_{ij}^2)^{-1/2}.\] By integrating this yields the formula
\[\tau_{ij} = \frac{1}{2\pi} \arcsin \frac{\sigma_{ij}}{2}.\]
To connect back to Kruskal's result, since we already assumed $\sigma_{ii} = \sigma_{jj} = 1$, then the correlation $\rho_{ij} = \sigma_{ij}$. Now, if we instead take the transformation to be be uniform on $[-1,1]$, we would find 
$\tau_{ij} = \frac{1}{\pi} \arcsin \frac{\sigma_{ij}}{2} = \frac{1}{\pi} \arcsin \frac{\rho_{ij}}{2}$
and that the variance of each is $\tau_{ii} = 1/6.$ Thus the correlation is
\[ \rho'_{ij} = \frac{\tau_{ij}}{\sqrt{\tau_{ii} \tau_{jj}}} = \frac{6}{\pi} \arcsin \frac{\rho_{ij}}{2} .\]
\end{example}

\subsection{Fourier series}
We now assume that our functions $f_i$ are given by a Fourier series
$$f_i(x)\ = \sum_{-K}^{K}a^i_n\, e^{i n x}$$
where 
$$a^i_n = \frac{1}{2\pi} \int_{0}^{2 \pi} f_i(x)\, e^{-i n x}\,dx.$$ We now follow the same process of the previous subsection. We only sketch the details since the proof uses the same idea. (In fact the following theorem relies on the series method of Section~\ref{sec:comp.smooth}, but the proof resembles that of the the transform method, so we include it here.)
\begin{theorem}
Suppose the functions $f_i$ have Fourier series of the above form. Then 
$$\tau_{ij} = \sum_n \sum_m a_n^i\,a_m^j (e^{-nm\sigma_{ij}} -1)\, e^{-\frac{1}{2}(n^2 \sigma_{ii}^2 + m^2 \sigma_{jj}^2)}.$$
\end{theorem}
\proof
We only note that 
$$f_i^{(2u+k)}(0) = \sum_n a^i_n (-i)^{2u+k} n^{2u+k}$$ and then the argument is essentially the same as in the Fourier transform case.
\endproof 
\begin{example}
Let $f_i = f_j = \sin x + \sin 2x.$ Then $$a_{-2} = \frac{-1}{2 i}, \,\,\,a_{-1} = \frac{-1}{2 i}, \,\,\,a_{1} = \frac{1}{2 i},\,\,\, a_{2} = \frac{1}{2 i}.$$ Hence
$$\tau_{ij} = \sinh (4\, \sigma_{ij})e^{-2(\sigma_{ii}^2 + \sigma_{jj}^2)}+ \sinh (2 \,\sigma_{ij})e^{-(2\sigma_{ii}^2 + \sigma_{jj}^2/2)} +\sinh (2 \,\sigma_{ij})e^{-(\sigma_{ii}^2/2 + 2\sigma_{jj}^2)} +\sinh  \sigma_{ij}e^{-(\sigma_{ii}^2 + \sigma_{jj}^2)/2}.$$

\end{example}

\subsection{The Laplace transform}
\begin{theorem}
Let $$f_i(x) = \int_0^{\infty}g_i(t) e^{-xt} dt.$$ We assume that $g_i$ has compact support. Then there exists $C_i$ and $M_i$ such that 
$$\int_0^{\infty}|g_i(t)\,t^k | e^{-xt} dt \leq C_iM_i^k$$ for all $i$ and $k$. It follows that 
$$\tau_{ij} = \int_{0}^{\infty} \int_{0}^{\infty} g_i(t_1)\,g_j(t_2)\, (e^{-t_1\,t_2\,\sigma_{ij}}-1)\, e^{-(t_1\sigma_{ii} +t_2 \sigma_{jj})/2}\,dt_1dt_2.$$
\end{theorem}
We leave the proof to the reader as it is the same as in the previous two cases.
(The reason for the compact support condition to to guarantee that the Laplace transform is defined for all $x$.)

\begin{center}* * * \end{center}

To summarize so far, we have (1)~formulas for smooth functions whose derivatives satisfy (\ref{der.est}) that yields a nice power series expansion for $\tau_{ij}$ in terms of $\sigma_{ij}$ and (2)~an integral representation for functions that are Fourier transforms with a series expression whose coefficients are made up of integrals. In the latter case, the transform functions may be discontinuous.

If $f_i$ is a polynomial, then the series method is very convenient. If $f_i$ is the Fourier or Laplace transform of some known function, then the transform method is convenient. If $f_i(x) = \sin x/x,$ either method can be used since in this case, the derivatives at zero are bounded by $1$ and $f_i$ is the Fourier transform of a characteristic function. 

Before we end this section, we should point out that (\ref{tau.int}) tells us that if 
$f_n$ is a sequence of functions that approaches $f$ in any reasonable sense (pointwise, $L^1$ sense) then the corresponding 
$(t_{ij})_n$ approach $\tau_{ij}.$ Thus one can always use an approximation that satisfies the conditions of one of our methods to approximate $\tau_{ij}.$

\section{Estimating the covariance entries}\label{sec:est}
In certain applications it is important to have estimates for $\tau_{ij}$ in terms of the $\sigma_{ij}.$ See for example \cite{baptista2024learning}.
We will consider these in two different ways, one using the techniques of the power series approach and the other using the Fourier transform.
\begin{theorem}
Let $f = f_i = f_j$ satisfy the condition \ref{der.est} with constants $C$ and $M.$
Then
\[\tau_{ij} = G_{1ij}(1/2) \sigma_{ij} + \frac{1}{2}\,G_{2ij}(1/2) \sigma_{ij}^2 +O(\sigma_{ij}^3).\] Here the constant for the $O$ estimate depends only on $f, \sigma_{ii},$ and $ \sigma_{jj}.$

\end{theorem}
\proof
We have \[\tau_{ij} - G_{1ij}(1/2) \sigma_{ij} -\frac{1}{2}\,G_{2ij}(1/2) \sigma_{ij}^2 = \sum_{k \geq 3} \frac{G_{kij}(1/2)\sigma_{ij}^k}{k!}.\]
Now $F_{kij}(x) = \sum_u \frac{f^{(2u +k)}(0) x^u}{u!}$ and is thus bounded by
$C M^k e^{M^2 x}.$ This yields a bound of 
\[ C^2 e^{M^2 (\sigma_{ii }+ \sigma_{jj})/2} \sum_{k\geq 3}\frac{M^{2k}\,|\sigma_{ij}|^k}{k!}\] for the right-hand side of the above displayed equation. Since $|\sigma_{ij}|$ is bounded by $N = \max \{\sigma_{ii}, \sigma_{jj}\}$ we can conclude that for our error term we have a bound of 
\[\frac{1}{6}\sigma_{ij}^3\, C^2\, M^6\, e^{M^2 (\sigma_{ii }+ \sigma_{jj})/2}\, e^{N M^2}.\]
\endproof
If we know more about $\sigma_{ij},$ say that it is at most $\epsilon,$ then the above bound can be improved to 
\[\frac{1}{6}\sigma_{ij}^3 \,C^2\, M^6 \,e^{M^2 \,(\sigma_{ii }+ \sigma_{jj})/2} \,e^{\epsilon M^2}.\]

\begin{theorem}
Let $f = f_i = f_j, i\neq j$ be in $L^2(-\infty, \infty)$ with inverse transform $g.$ Suppose $g$ is bounded by a constant $N.$ Also suppose that $$\frac{|\sigma_{ij}|}{\sqrt{\sigma_{ii}}\sqrt{\sigma_{jj}}} < a <1.$$
Then
\[\tau_{ij} = G^*_{1ij}(1/2) \sigma_{ij} + \frac{1}{2}\,G^*_{2ij}(1/2) \sigma_{ij}^2 + O(\sigma_{ij}^3),\]
where the $O$ depends on $f, \sigma_{ii}, \sigma_{jj}$ and $a.$

(Recall $G^*_{kij} = F^*_{ki}(\sigma_{ii}x)F^*_{kj}(\sigma_{jj}x)$ where $F^*_{kj}(x) = \frac{(-i)^k}{2\pi} \int_{-\infty}^{\infty} g_j(y)\,\,y^k e^{-y^2 x} dy.$)
\proof
We have that 
\[ \tau_{ij}  = \frac{1}{4\pi^2} \int_{-\infty}^\infty \int_{-\infty}^\infty g(y)\,g(z)\,(e^{- y\,z\,\sigma_{ij}}-1)\, e^{(-y^2\sigma_{ii} -z^2 \sigma_{jj})/2}\,dydz.\]
Expanding the middle exponential yields
$$\tau_{ij} - G^*_{1ij}(1/2) \sigma_{ij} -  \frac{1}{2}\,G^*_{2ij}(1/2) \sigma_{ij}^2= \frac{1}{4\pi^2}\sum_{k\geq 3} \frac{(-\sigma_{ij})^k}{k!}\int_{-\infty}^\infty \int_{-\infty}^\infty g(y)\,g(z)\,y^k z^k\, e^{(-y^2\sigma_{ii} -z^2 \sigma_{jj})/2}\,dydz.$$
Using the bounds we found for the integrals
\[\int_{-\infty}^\infty  g(y)\,y^k \, e^{(-y^2\sigma_{ii}/2 )}\,dy\] in (\ref{tau.int}) it follows that the $O$ is at most $\frac{N^2}{\pi^2 \sqrt{\sigma_{ii}}\sqrt{\sigma_{jj}}(1 -a)}.$

\endproof

\end{theorem}

\begin{example}

Let $f_i(x) = f_j(x) = f(x) = \frac{1}{1+x^2}.$ This function does not satisfy  condition~(\ref{der.est}) since the convergence of the power series has radius $1$. 
However $f(x)$ is the Fourier transform of $\pi e^{-|y|}$ and thus we can write 
\[ \tau_{ij}  = \frac{1}{4} \int_{-\infty}^\infty \int_{-\infty}^\infty e^{-|y|}e^{-|z|}\,(e^{- y\,z\,\sigma_{ij}}-1)\, e^{(-y^2\sigma_{ii} -z^2 \sigma_{jj})/2}\,dydz.\]

We wish to compute $G^*_{1ij}(1/2) \sigma_{ij} + \frac{1}{2}\,G^*_{2ij}(1/2) \sigma_{ij}^2.$ Since our function is even $G^*_{1ij}(1/2) = 0.$
We have
\[ F^*_{2}(x) = -\int_{0}^\infty e^{-y}\,\,y^{2} e^{-y^2 x} dy.\]
Using integration by parts we find this yields that
\[ G^*_{2ij} = A(\sigma_{ii})A(\sigma_{jj})\]
where
\[ A(z) = \frac{1}{z} \left(  1 - \frac{\pi (1 +z)}{\sqrt{2 z}} e^{1/2} [1 - \Phi(1/\sqrt{2z})]  \right)\]
and $\Phi(z) = \frac{2}{\sqrt{\pi}}\int_0^z e^{-t^2} dt$ is the error function.
\end{example}

\section{Glossary}\label{sec:gloss}
Here is a list of basic properties of transformed variables:
\begin{enumerate}
\item $$ \nu_i = \frac{1}{\sqrt{2\,\pi \,\sigma_{ii}}}\int_{-\infty}^\infty  f_i(y)\, e^{(-y^2/2\sigma_{ii})}\,dy,$$
    \item For $i\neq j,$ $$\tau_{ij} =
     c \int_{-\infty}^{\infty}\int_{-\infty}^{\infty} f_i(x)f_j(y) e^{-\frac{\sigma_{jj}x^2}{2 d}}e^{\frac{\sigma_{ij}xy}{ d}}e^{-\frac{\sigma_{ii}y^2}{2 d}}dx\, dy - \nu_i\nu_j,\,\,\,
 $$ $$d = \sigma_{ii}\sigma_{jj} - \sigma_{ij}^2 \,\,\,(d >0), \,\,c = \frac{1}{d^{1/2}2\pi}.$$
    \item For $i=j,$ $$\tau_{ij} = \frac{1}{\sqrt{2\,\pi \,\sigma_{ii}}}\int_{-\infty}^{\infty} f_i(x)^2 e^{(-x^2/2\sigma_{ii})} dx -\nu_i^2  $$
\end{enumerate}
Define 
$$\tau_{ij}(f_i, f_j) =
     c \int_{-\infty}^{\infty}\int_{-\infty}^{\infty} f_i(x)f_j(y) e^{-\frac{\sigma_{jj}x^2}{2 d}}e^{\frac{\sigma_{ij}xy}{ d}}e^{-\frac{\sigma_{ii}y^2}{2 d}}dx\, dy - \nu_i\nu_j,\,\,\,$$ for $i\neq j$ and for $i =j $ by 
     $$\tau_{ij}(f_i,f_i) = \frac{1}{\sqrt{2\,\pi \,\sigma_{ii}}}\int_{-\infty}^{\infty} f_i(x)^2 e^{(-x^2/2\sigma_{ii})} dx -\nu_i^2  .$$
  Then we have the additional properties
  \begin{enumerate} \setcounter{enumi}{3}
\item $\tau_{ij}(f_i +g_i, f_j) = \tau_{ij}(f_i, f_j) + \tau_{ij}(g_i, f_j)$
\item $\tau_{ij}(f_i, f_j + g_j) = \tau_{ij}(f_i, f_j) + \tau_{ij}(f_i, g_j)$
\item $\tau_{ij}(c_if_i, c_jf_j) = c_ic_j \tau_{ij}(f_i, f_j)$ assuming $c_i$ and $c_j$ are constants
\item Let $ f_i = f_{io} + f_{ie}\,\,\,\,f_j = f_{jo} + f_{je}$ where $f_{io},\,\, f_{ie}$ are odd and even functions. Then
$\tau_{ij}(f_i, f_j) = \tau_{ij}(f_{io}, f_{jo}) + \tau_{ij}(f_{ie}, f_{je}).$
\item 
Let $$f_k(x) = \hat{g}_k(x) = \frac{1}{2\pi}\int_{-\infty}^\infty g_k(y) e^{-i\,x\,y}\,dy,$$ for $k = i,j.$ Then
$$\tau_{ij} = \frac{1}{4\pi^2} \int_{-\infty}^\infty \int_{-\infty}^\infty g_i(y)\,g_j(z)\,(e^{- y\,z\,\sigma_{ij}}-1)\, e^{(-y^2\sigma_{ii} -z^2 \sigma_{jj})/2}\,dydz.$$
  \end{enumerate}
To end this section, here is a table of some known transformations when $f_i = f_j$. This table does not include all the examples in the paper, and it does contain some new ones. If the example is in the earlier text, the example number is given. The $\tau_{ii}$ column is only filled if a different expression is needed than for $\tau_{ij}$; otherwise, one can assume that the $\tau_{ij}$ formula still holds for $i = j$. The first ten examples were computed with the series method; the last four (starting with the Gaussian function) rely on the transform method.

\renewcommand{\arraystretch}{1.8}
\begin{table}[H]
\hspace{-2em}
 \begin{tabular}{c|c|c|c|c}
& $f_i$ & $\tau_{ij}$ & $\tau_{ii}$ & $\text{Example in text}$\\\hline 
1& $\sin x$  & $e^{-\frac{\sigma_{ii}+\sigma_{jj}}{2}}\sinh \sigma_{ij} $& - & \\\hline 
2&$\cos x$ & $e^{-\frac{\sigma_{ii}+\sigma_{jj}}{2}}(\cosh \sigma_{ij} -1)$ & - &\ref{ex:cos} \\\hline 
3&$\sin ax$ & $e^{-a^2\frac{\sigma_{ii}+\sigma_{jj}}{2}}\sinh a^2 \sigma_{ij}$ & -& \\\hline 
4&$\cos ax$ & $e^{-a^2\frac{\sigma_{ii}+\sigma_{jj}}{2}}(\cosh a^2\sigma_{ij} -1)$ &-&   \\\hline 
5&$\sinh ax$ &  $e^{a^2\frac{\sigma_{ii}+\sigma_{jj}}{2}}\sinh a^2\sigma_{ij}$ & -& \\\hline 
6&$\cosh ax$ & $e^{a^2\frac{\sigma_{ii}+\sigma_{jj}}{2}}(\cosh a^2\sigma_{ij} -1)$& - &  \\\hline 
7&$e^x$ & $e^{\frac{\sigma_{ii}+\sigma_{jj}}{2}}(e^{\sigma_{ij}} -1)$ & -& \\\hline 
8&$e^{ax}$ & $e^{a^2\frac{\sigma_{ii}+\sigma_{jj}}{2}}(e^{a^2\sigma_{ij}} -1)$& -& \\\hline 
9&$\frac{x^{2n}}{(2n)!}$&$\frac{\sigma_{ii}^n\sigma_{jj}^n }{2^{2n}}\sum_{k=1}^n \frac{1}{(n-k)!^2 (2k)!}(\frac{4\sigma_{ij}^2}{\sigma_{ii}\sigma_{jj}})^k$ & -& \\\hline 
10&$\frac{x^{2n+1}}{(2n+1)!}$ & $\frac{\sigma_{ii}^n\sigma_{jj}^n \sigma_{ij}}{2^{2n}}\sum_{k=0}^n \frac{1}{(n-k)!^2 (2k+1)!}(\frac{4\sigma_{ij}^2}{\sigma_{ii}\sigma_{jj}})^k$ & -& \\\hline 
11&$\frac{1}{2 \pi a} e^{-x^2 /(2a)}$& $\frac{1}{2\pi}(((a+\sigma_{ii})(b+\sigma_{jj}) - \sigma_{ij}^2)^{-1/2}$ & - &\ref{ex:gauss}\\ 
&&  $-((a+\sigma_{ii})(b+\sigma_{jj}))^{-1/2})$ &  \\\hline
12&$\chi_{[-1,1]} (\sigma_{ii} = \sigma_{ij} =1)$ & $\frac{1}{\pi \, e}\sum_{k = 1}^{\infty}\,\frac{\sigma_{ij}^{2k}}{(2k)!\,\,2^{2k-2}}\,(H_{2k-1}(1/\sqrt{2}))^2$ &  $\Phi(1/\sqrt{2}) - \Phi(1/\sqrt{2})^2$ &\ref{ex:char}\\\hline
13&$x \cdot \chi_{[-1,1]} (\sigma_{ii} = \sigma_{ij} =1) $ &$\frac{1}{\pi^2}\left(\pi\Phi(1/\sqrt{2}) - \left(\frac{2\pi}{e}\right)^{1/2}\right)^2\sigma_{ij} $ & $\Phi(1/\sqrt{2}) - \sqrt{\frac{2}{\pi e}}$&\ref{ex:char-x}\\
&&$+ \frac{2}{ \pi e} \sum_{k=2}^{\infty} \frac{\sigma_{ij}^{2k+1}(2k-1)^2}{(2k+1)!2^{2k-2}}H_{2k-2}(1/\sqrt{2})^2$ & &\\
\hline
14&$\frac{1}{\sqrt{2 \pi}} \int_{-\infty}^x e^{-t^2/2}$ & $\frac{1}{2\pi}\arcsin{\frac{\sigma_{ij}}{2}}$ & - &\ref{ex:uni}\\
\end{tabular}\caption{Theoretical results for some common transform functions.\label{tab:t}}
\end{table}

\section{Numerical examples}\label{sec:num}
Here are several numerical examples; all start with $\sigma_{ii} = 1$ and $\sigma_{ij} = 1/4$ for $i \neq j$. Both tables~\ref{tab:tvn} and~\ref{tab:bad} define a specific function (with parameters fixed), and symbolic expressions for $\tau_{ij}$ as well as $\tau_{ii}$ if needed. These expressions are then evaluated numerically, and also compared to empirical estimates, given $10^6$ samples. In table~\ref{tab:tvn}, there is overall very close agreement between the theoretical and empirical answers, usually up to the third decimal place, or more. For some functions (or choices of parameters), however, the function perturbs the data more significantly, and the empirical estimates may suffer, as seen in table~\ref{tab:bad}. Although these estimates would probably improve with standardization of the samples, the theoretical answers are still reliable and do not depend on standardization or sample size.

The Python code used for the numerical results is available here~\cite{morrison2024nonpararnormal}.

\begin{table}[H]
\hspace{-4em}
 \begin{tabular}{c|c|c|c|c|c|c|c}
&$f_i$ & $\tau_{ij}$ & $\tau_{ij}$& $\tau_{ij}$ & $\tau_{ii}$& $\tau_{ii}$& $\tau_{ii}$\\
&&  &{\footnotesize(evaluated)}& {\footnotesize(empirical)}& & \text{{\footnotesize(evaluated)}}& \text{{\footnotesize(empirical)}} \\\hline 
1&$\sin x$  & $e^{-\frac{\sigma_{ii}+\sigma_{jj}}{2}}\sinh \sigma_{ij}$ & 0.0929& 0.0916& - & 0.4323& 0.4307 \\\hline 
2&$\cos x$ & $e^{-\frac{\sigma_{ii}+\sigma_{jj}}{2}}(\cosh \sigma_{ij} -1) $& 0.0116&0.0118& - & 0.1998& 0.2001 \\\hline 
3&$\sin ax, \,\, a = 2 $& $e^{-a^2\frac{\sigma_{ii}+\sigma_{jj}}{2}}\sinh a^2 \sigma_{ij}$ & 0.0215& 0.0228&-  &0.4998& 0.5004 \\\hline 
4&$\cos ax, \,\,a = 1/2$ & $e^{-a^2\frac{\sigma_{ii}+\sigma_{jj}}{2}}(\cosh a^2\sigma_{ij} -1)$ & 0.0015 & 0.0017& - & 0.0245& 0.0244  \\\hline 
5&$\sinh ax,\,\,a=1 $& $ e^{a^2\frac{\sigma_{ii}+\sigma_{jj}}{2}}\sinh a^2\sigma_{ij}$ & 0.6867& 0.7066& - & 3.1945& 3.2154 \\\hline 
6&$\cosh ax,\,\,a=3/2$ & $e^{a^2\frac{\sigma_{ii}+\sigma_{jj}}{2}}(\cosh a^2\sigma_{ij} -1)$&1.5410& 1.3411& - & 36.0208& 36.8717  \\\hline 
7&$e^x$ & $e^{\frac{\sigma_{ii}+\sigma_{jj}}{2}}(e^{\sigma_{ij}} -1)$ & 0.7721&0.7523&- & 4.6708& 4.6282\\\hline 
8&$e^{ax},\,\,a=1/3$ & $e^{a^2\frac{\sigma_{ii}+\sigma_{jj}}{2}}(e^{a^2\sigma_{ij}} -1)$& 0.0315& 0.0323& - & 0.1313& 0.1317 \\\hline 
9&$\frac{x^{2n}}{(2n)!},\,\, n=1$&$\frac{\sigma_{ii}^n\sigma_{jj}^n }{2^{2n}}\sum_{k=1}^n \frac{1}{(n-k)!^2 (2k)!}(\frac{4\sigma_{ij}^2}{\sigma_{ii}\sigma_{jj}})^k$ & 0.0312& 0.0302& - & 0.5000& 0.5019 \\\hline 
10&$\frac{x^{2n+1}}{(2n+1)!},\,\,n=2$ & $\frac{\sigma_{ii}^n\sigma_{jj}^n \sigma_{ij}}{2^{2n}}\sum_{k=0}^n \frac{1}{(n-k)!^2 (2k+1)!}(\frac{4\sigma_{ij}^2}{\sigma_{ii}\sigma_{jj}})^k$ & 0.0046& 0.0046& - &0.0656& 0.0690 \\\hline 
11&$\frac{1}{2 \pi a} e^{-x^2 /(2a)},\,\,a=1$& $\frac{1}{2\pi}(((a+\sigma_{ii})(b+\sigma_{jj}) - \sigma_{ij}^2)^{-1/2}$ & 0.0006 & 0.0006 & - & 0.0123& 0.0123  \\
&&  $-((a+\sigma_{ii})(b+\sigma_{jj}))^{-1/2})$ &&&&&  \\\hline
12&$\chi_{[-1,1]}$,  & $\frac{1}{\pi \, e}\sum_{k = 1}^{\infty}\,\frac{\sigma_{ij}^{2k}}{(2k)!\,\,2^{2k-2}}\,(H_{2k-1}(1/\sqrt{2}))^2$ & 0.0075 & 0.0073& $\Phi(1/\sqrt{2})$ & 0.2166 & 0.2160 \\
&$\sigma_{ii} = \sigma_{jj} =1$ &&&& $- \Phi(1/\sqrt{2})^2$   && \\ \hline
13&$x \cdot \chi_{[-1,1]}$, &$\left(\Phi(1/\sqrt{2}) - \sqrt{\frac{2}{\pi e}}\right)^2\sigma_{ij}$  & 0.0099 & 0.0092& $\Phi(1/\sqrt{2})$ &0.1987& 0.1973 \\
&$\sigma_{ii} = \sigma_{ij} =1$  &$+ \frac{2}{ \pi e} \sum_{k=2}^{\infty} \frac{\sigma_{ij}^{2k+1}(2k-1)^2}{(2k+1)!2^{2k-2}}H_{2k-2}(1/\sqrt{2})^2$ & && $- \sqrt{\frac{2}{\pi e}}$&&\\ \hline
14&$\frac{1}{\sqrt{2 \pi}} \int_{-\infty}^x e^{-t^2/2}$ & $\frac{1}{2\pi}\arcsin{\frac{\sigma_{ij}}{2}}$ & 0.0199& 0.0194&-&0.0833&0.0831 \\
\end{tabular}
\begin{center}
    \caption{Theoretical answers evaluated for some specific functions, and compared to numerical estimates based on $10^6$ samples. There is very good agreement between the two.\label{tab:tvn}}\end{center}
\end{table}

\begin{table}[H]
\centering
 \begin{tabular}{c|c|c|c|c|c|c|c}
&$f_i$ & $\tau_{ij}$ & $\tau_{ij}$& $\tau_{ij}$ & $\tau_{ii}$ & $\tau_{ii}$& $\tau_{ii}$\\
&&  &{\footnotesize(evaluated)}& {\footnotesize(empirical)}& & {\footnotesize(evaluated)}& \text{{\footnotesize(empirical)}} \\\hline 
3&$\sin ax, \,\, a =3$ & $e^{-a^2\frac{\sigma_{ii}+\sigma_{jj}}{2}}\sinh a^2 \sigma_{ij}$ & 0.0006& -0.0018&-  &0.5000& 0.4978 \\\hline 
6&$\cosh ax,\,\,a=5/2$ & $e^{a^2\frac{\sigma_{ii}+\sigma_{jj}}{2}}(\cosh a^2\sigma_{ij} -1)$&771.9& 164.2& - & 133651& 5311965  \\\hline 
8&$e^{ax},\,\,a=2$ & $e^{a^2\frac{\sigma_{ii}+\sigma_{jj}}{2}}(e^{a^2\sigma_{ij}} -1)$& 93.82& 66.97& - & 2926& 2158 \\
\end{tabular}\caption{Numerical estimates diverge from theory for some functions.\label{tab:bad}}
\end{table}

\section{Discussion}\label{sec:conc}
Diagonal (or marginal) transformations of multivariate Gaussian variables appear in diverse fields including probabilistic graphical models~\cite{harris2013pc, liu2012nonparanormal}, information theory~\cite{singh2017nonparanormal}, causal effect estimation~\cite{mahdi2018estimating}, and inference~\cite{mulgrave2018bayesian}. In this paper, we show two ways to compute first and second moments of these distributions exactly, for a broad class of transformation functions. We also provide some estimates that bound the covariance between the transformed variables in terms of the original multivariate normal covariance.

In addition to the theoretical results, we provide a glossary of identities about these computations, a table of common transformation functions and the associated moments, and numerical results that compare our exact answers to empirical estimates. When the functions are fairly well-behaved, there is excellent agreement between our theoretical and empirical results. On the other hand, certain functions perturb the data such that numerical computations of the covariance matrix suffer from high variance; in these cases, the theoretical results give much more reliable answers. And of course in all experiments, the exact answers do not depend on any sample size effects.

What about higher moments? In principle, these methods should readily extend to compute those as well. In the series approach, this would require counting higher mixed moments of a multivariate normal (e.g., with Wick's theorem), which is not an easy combinatorial problem. Perhaps more accessible would be higher moments via the transform method, which would instead require higher-dimensional integration. Ideally one would then summarize all of the moments for the transformed random variable through some type of moment generating function, yielding not just a collection of summary statistics of the new distribution but rather a complete analytic description.

\section*{Acknowledgments}The authors would like to thank Elmar Plischke for pointing out the connection between our work and the papers by Kruskal, and Clemen and Reilly. The authors also thank Ricardo Baptista for several helpful comments and suggestions. RM acknowledges the Johnson \& Johnson Foundation and the Women in STEM2D Awards Program.

\bibliographystyle{plain}
\bibliography{references.bib}


\end{document}